\documentclass[10pt]{article}
\usepackage{amsmath}
\usepackage{latexsym}
\usepackage{amssymb}
\usepackage{amsfonts}
\usepackage{amsthm}
\title{GENERAL RELATIVITY AND GAUGE THEORY: \\BEYOND THE MIRROR}
\author{J.-F. Pommaret \\  jean-francois.pommaret@wanadoo.fr  \\
ORCID:0000-0003-0907-2601 \\ }
\date{  }
\begin{document}
\maketitle

\thispagestyle{empty}
\noindent
{\bf ABSTRACT}\\

Lie pseudogroups are groups of transformations solutions of systems of ordinary (OD) or partial differential (PD) equations. The purpose of this paper is to present an elementary summary of a few recent results obtained through the application of the formal theory of systems of OD or PD equations and Lie pseudogroups to engineering (elasticity, electromagnetism) or mathematical physics (general relativity, gauge theory) and their couplings (piezoelectricity, photoelasticity). The work of Cartan is superseded by the use of the canonical Spencer sequence while the work of Vessiot is superseded by the use of the canonical Janet sequence but the link between these two sequences and thus these two works is still not known today. Using differential duality in the linear framework, the adjoint of the Spencer operator for the group of conformal transformations provides the Cosserat equations, the Maxwell equations and the Weyl equations on equal footing. Such a result allows to unify the finite elements of engineering sciences but also leads to deep contradictions in the case of gravitational waves. Indeed, the Beltrami operator (1892) which is parametrizing the Cauchy operator of elasticity by means of 6 stress functions is nothing else than the self-adjoint Einstein operator (1915) in dimension 3 for the deformation of the metric which is parametrizing the {\it div} operator induced from the Bianchi identities. The same confusion between the Cauchy and {\it div} operators is existing on space-time as the Cauchy operator can be parametrized by the adjoint of the Ricci operator. Accordingly, the foundations of engineering and mathematical physics must be revisited within this new framework, though striking it may sometimes look like.

\vspace{3cm}

Differential duality, Maxwell equations, Einstein equations, Gravitational waves, 

\vspace{4cm}
\newpage

\noindent
{\bf 1) INTRODUCTION}  \\

The purpose of this self-contained but difficult paper is to revisit {\it general relativity} (GR) and {\it gauge theory} (GT) in view of the latest mathematical developments existing today in {\it group theory}, {\it system theory} and {\it module theory}, namely:  \\ 

\noindent
$\bullet$ {\it Systems} : The order of the successive operators appearing in the conformal Killing resolution highly depend on the dimension $n$, a result recently confirmed by A. Quadrat (INRIA) while using computer algebra ([21]). They are respectively $ 3 \underset{1}{\longrightarrow} 5 \underset{3}{\longrightarrow }5 \underset{1}{ \longrightarrow }3 \rightarrow 0 $ when $n=3$,  
  $  4 \underset{1}{\longrightarrow} 9 \underset{2}{\longrightarrow} 10 \underset{2}{\longrightarrow} 9 \underset{1}{\longrightarrow} 4  \rightarrow 0  $ when $n=4$ and finally                                         
  $   5 \underset{1}{\longrightarrow} 14 \underset{2}{\longrightarrow} 35 \underset{1}{\longrightarrow} 35 \underset{2}{\longrightarrow} 14 \underset{1}{\longrightarrow} 5 \rightarrow 0  $ when $n=5$. This result leads to revisit conformal geometry.  \\

\noindent
$\bullet$ {\it Groups} : The Ricci and the Maxwell tensors have only to do with the second order jets of the conformal group, called {\it elations} by Cartan (1922) ([18, 23, 25, 26]). This result questions the mathematical foundations of both general relativity (GR) and electromagnetism (EM]).  \\

\noindent
$\bullet$ {\it Modules} : Contrary to the Maxwell equations, the Einstein equations cannot be parametrized themselves. Such a result is not coherent with the use of " {\it extension modules} " in homological algebra. As a byproduct, the Cauchy stress equations must not be confused with the divergence-type condition for the Einstein tensor obtained by contracting the Bianchi identities ([13, 15, 19, 22, 23, 33]). This result questions the origin and existence of gravitational waves ([28]). \\

We briefly recall the historical framework leading to these new results (Compare to [19]).  \\
The concept of "{\it group}" has been introduced in mathematics for the first time by E. Galois (1830) and slowly passed from algebra to geometry with the work of  S. Lie on {\it Lie groups} (1880) and {\it Lie pseudogroups} (1890) of transformations. The concept of a finite length {\it differential sequence}, now called {\it Janet sequence}, has been described for the first time as a footnote by M. Janet (1920). Then, the work of D. C. Spencer (1970) has been the first attempt to use the formal theory of systems of partial differential equations in order to study the formal theory of Lie pseudogroups. However, the linear and nonlinear {\it Spencer sequences} for Lie pseudogroups, though never used in physics, largely supersede the "{\it Cartan structure equations} " (1905) and are quite different from the "{\it Vessiot structure equations} " (1903), introduced for the same purpose but still not known today because they have never been acknowledged by E. Cartan or successors ([5, 27, 31]).   \\
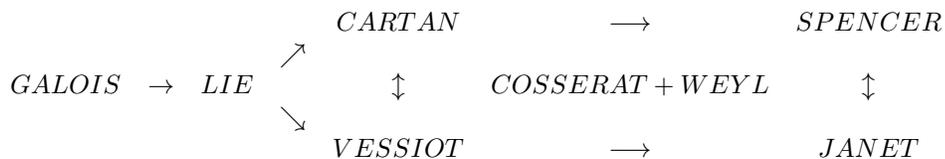
\[   \begin{array}{rcccccc}
                                         &   &  &   &  CARTAN   &   \longrightarrow   &  SPENCER  \\
                                  &    &     &   \nearrow  &   &   &   \\
      GALOIS & \rightarrow   & LIE   &   &   \updownarrow   & {COSSERAT}+{WEYL} &   \updownarrow   \\
                                  &   &   &   \searrow  &  &  &   \\
                                         &  &  &   &   VESSIOT   &   \longrightarrow   &  JANET   
   \end{array}     \]

Meanwhile, mixing differential geometry with homological algebra, M. Kashiwara (1970) has created "{\it differential homological algebra} ",  in order to study {\it differential modules} by means of {\it double duality} and the corresponding {\it extension modules} (See [14] for references and Zbl 1079.93001). \\

{\it By chance}, unexpected arguments have been introduced by the brothers E. and F. Cosserat (1909) in order to revisit elasticity and by H. Weyl (1918) in order to revisit electromagnetism through a {\it unique differential sequence} only depending on the structure of the {\it conformal group}. However, while the Cosserat brothers were only using ({\it translations} + {\it rotations}), Weyl has only been dealing with ({\it dilatation} + {\it elations}) as we shall explain ([2, 11, 32]). \\

The initial motivation for studying the methods used in this paper has been a $1000 \$ $ challenge proposed in $1970$ by J. Wheeler in the physics department of Princeton University while the author of this paper was a visiting student of D.C. Spencer in the close-by mathematics department:  \\ 
{\it Is it possible to express the generic solutions of Einstein equations in vacuum by means of the derivatives of a certain number of arbitrary functions, like the potentials for Maxwell equations ?}.\\
After recalling the negative answer we already provided in 1995 ([13]), the main purpose of this paper is to use the new techniques of {\it differential double duality} in order to revisit the mathematical foundations of general relativity and gauge theory that are leading to gravitational waves ([28]). We point out the fact that all the formulas presented could be obtained by computer algebra while using recent packages developed for this purpose by A. Quadrat and collaborators.\\

\noindent
{\bf 2) PARAMETRIZATION}  \\

Starting with the well known linear map $C : S_2T^* \rightarrow S_2T^*: R_{ij} \rightarrow E_{ij}=R_{ij}- \frac{1}{2}{\omega}_{ij}tr(R)$ between symmetric covariant tensors, where $\omega$ is a metric with $det(\omega)\neq 0$ and $tr(R)={\omega}^{rs}R_{rs}$, we may introduce the linear second order operators $Ricci: \Omega \rightarrow R$ and $Einstein: \Omega \rightarrow E$ obtained by linearization over $\omega$ and we have the relation $Einstein = C \circ Ricci $ where $C$ does not depend on any conformal factor ([3]).We recall the method used in any textbook for studying gravitational waves, which "{\it surprisingly} " brings the same map $ C :\Omega \rightarrow \bar{\Omega}=\Omega -\frac{1}{2}\omega \,tr(\Omega)$ in order to introduce the key unknown composite operator ${\cal{X}}:\bar{\Omega}\rightarrow \Omega \rightarrow E$, having therefore 
$Einstein = {\cal{X}} \circ C$. \\

The first goal will be to prove that only the use of {\it differential homological algebra}, a mixture of differential geometry (differential sequences, formal adjoint) and homological algebra (module theory, double duality, extension modules) {\it totally unknown by physicists}, is able to explain why the Einstein operator (with 6 terms) defined above is useless as it can be replaced by the Ricci operator (with 4 terms) in the search for gravitational waves equations. Indeed, using the fact that the {\it Einstein} operator is self-adjoint ([19]) when $\omega$ is the Minkowski metric (contrary to the {\it Ricci} operator as we shall see) and taking the respective (formal) adjoint operators (multiplying by test functions and integrating by parts), we get:  
\[  ad(Einstein)=ad(C) \circ ad({\cal{X}}) \Rightarrow Einstein= C \circ ad({\cal{X}})\Rightarrow ad({\cal{X}})=Ricci \Rightarrow 
{\cal{X}}= ad(Ricci)  \]
Meanwhile, the {\it Riemann} operator can be considered as an operator describing the (second order) {\it compatibility conditions} (CC) for the {\it Killing} operator $\xi \in T \rightarrow {\cal{L}}(\xi) \omega=\Omega \in S_2T^*$ with standard notations where ${\cal{L}}$ is the Lie derivative ([9, 27]). In this new framework, we shall prove that we no longer need to use the {\it Bianchi} operator as the first order CC for the {\it Riemann} operator. Also, we shall prove that the {\it relative parametrization} with $div$-type {\it differential constraints} needed in order to keep only the $Dalembert$ operator in the wave equations has {\it nothing to do} with any gauge transformation in the corresponding adjoint differential sequence, but has {\it only to do} with the search for a {\it minimal parametrization}, exactly like Maxwell in 1870 for elasticity ([20, 28]).\\

\noindent
{\bf EXAMPLE 2.1}: When $n=2$, the stress equations become ${\partial}_1{\sigma}^{11}+{\partial}_2{\sigma}^{12}=0, {\partial}_1{\sigma}^{21}+{\partial}_2{\sigma}^{22}=0$. Their second order parametrization ${\sigma}^{11}={\partial}_{22}\phi, {\sigma}^{12}={\sigma}^{21}=-{\partial}_{12}\phi, {\sigma}^{22}={\partial}_{11}\phi$ has been provided by George Biddell Airy (1801-1892) in 1863. It can be simply recovered in the following manner: \\
\[ \begin{array}{rcl}
{\partial}_1{\sigma}^{11}- {\partial}_2( - {\sigma}^{12})= 0 \hspace{5mm} & \Rightarrow & \hspace{5mm} \exists \varphi,\, {\sigma}^{11}={\partial}_2\varphi, {\sigma}^{12}= 
- {\partial}_1\varphi \\  
{\partial}_2{\sigma}^{22}- {\partial}_1( - {\sigma}^{21})=0  \hspace{5mm} & \Rightarrow & \hspace{5mm}  \exists \psi, \,{\sigma}^{22}={\partial}_1\psi, {\sigma}^{21}= - {\partial}_2\psi \\
   {\sigma}^{12}={\sigma}^{21} \Rightarrow {\partial}_1 \varphi - {\partial}_2\psi =0 \hspace{5mm}& \Rightarrow & \hspace{5mm}\exists \phi, \,\varphi={\partial}_2\phi, \psi={\partial}_1\phi  
   \end{array}   \]

When constructing a long prismatic dam with concrete, we may transform a problem of $3$-dimensional elasticity into a problem of $2$-dimensional elasticity by supposing that the axis $x^3$ of the dam is perpendicular to the river with ${\Omega}_{ij}(x^1,x^2), \forall i,j=1,2$ and ${\Omega}_{33}=0$ because the rocky banks of the river are supposed to be fixed. We may introduce the two {\it Lam\'{e} constants} $(\lambda,\mu)$ in order to describe the usual constitutive relations of an homogeneous isotropic medium as follows, passing from the standard case $n=3$ to the restricted case $ n=2$ just by setting:  \\
\[ {\sigma}=\frac{1}{2}\lambda \, tr(\Omega)\, {\omega} + \mu \, {\Omega}, \, \, tr(\Omega)={\Omega}_{11}+{\Omega}_{22} \hspace{4mm}  \Rightarrow   \hspace{4mm}  
 \mu \, {\Omega}=    \sigma - \frac{\lambda}{2( \lambda + \mu)}\, tr(\sigma) \, \omega , \, \, tr(\sigma)={\sigma}^{11} + {\sigma}^{22} \]
even though ${\sigma}^{33}=\frac{1}{2}\lambda ({\Omega}_{11}+{\Omega}_{22})=\frac{1}{2}\lambda tr(\Omega) $. Let us consider the {\it right square} of the diagram below with locally exact rows, where any vector bundle is simply denoted by its fiber dimension:  \\
\[   \begin{array}{rcccccl}
  & 2  & \stackrel{Killing}{\longrightarrow} & 3 & \stackrel{Riemann}{\longrightarrow} & 1  & \rightarrow 0 \\
  & \vdots &  & {\downarrow\uparrow} &  &  \vdots & \\
 0 \leftarrow &  2 & \stackrel{Cauchy}{\longleftarrow} & 3 & \stackrel{Airy}{\longleftarrow} & 1&
\end{array}  \]
Taking into account the linearization of the only component of the Riemann tensor over the Euclidean metric $\omega$ when $n=2$ and substituting the Airy parametrization, we obtain:  \\
\[ tr(R)\equiv d_{11}{\Omega}_{22}+d_{22}{\Omega}_{11}-2d_{12}{\Omega}_{12}=0 \hspace{3mm} \Rightarrow 
\hspace{3mm} \mu \, tr(R)\equiv \frac{\lambda + 2 \mu}{2(\lambda +\mu)} \Delta \Delta \phi=0 \hspace{3mm} \Rightarrow 
\hspace{3mm} \Delta \Delta \phi=0  \]
where the linearized {\it scalar curvature} $tr(R)$ is allowing to define the {\it Riemann operator} in the previous diagram, namely the only {\it compatibility condition} (CC) of the Killing operator. It remains to exhibit an arbitrary homogeneous polynomial solution of degree $3$ and to determine its $4$ coefficients by the boundary pressure conditions on the upstream and downstream walls of the dam. The Airy potential $\phi$ has {\it nothing to do} with the perturbation $\Omega$ of the metric $\omega$ and the {\it Airy} operator is nothing else but the adjoint of the $Riemann$ operator, that is $Airy=ad(Riemann)$.  \\

\noindent
{\bf EXAMPLE 2.2}: When $n=3$, we may now use the {\it left square} of the following diagram with locally exact rows:  \\
\[   \begin{array}{cccccc}
   &  3  & \stackrel{Killing}{\longrightarrow} & 6 & \stackrel{Riemann}{\longrightarrow} & 6   \\
  & \vdots &  & {\downarrow\uparrow} &  &  \vdots  \\
 0 \leftarrow &  3 & \stackrel{Cauchy}{\longleftarrow} & 6 & \stackrel{Beltrami}{\longleftarrow} & 6
\end{array}  \]
where the self-adjoint operator $Beltrami=ad(Riemann)$ has ben introduced by E. Beltrami in 1892. We may substitute the $3$-dimensional constitutive relations with Lam\'{e} constants $(\lambda, \mu)$ in the Cauchy stress equations and get, when $\vec{f}=\vec{g}$ ({\it gravity}) is now the right member:   \\
\[  (\lambda + \mu) \vec{\nabla}(\vec{\nabla}.\vec{\xi})+\mu \Delta \vec{\xi}=\vec{f}\hspace{2mm}  \stackrel{\vec{\nabla}}{\Rightarrow} \hspace{2mm}  (\lambda + 2\mu)\Delta tr(\Omega) =0 \Rightarrow \Delta tr(\Omega)=0\Rightarrow \Delta tr(\sigma)=0    \]
We discover at once that the origin of elastic waves is shifted by {\it one step backwards}, {\it from the right square to the left square} of the diagram. Indeed, using inertial forces $\vec{f}=\rho \,{\partial}^2\vec{\xi}/\partial t^2$ for a medium with mass $\rho$ per unit volume in the right member of Cauchy stress equations because of Newton law and the vector identity $\vec{\nabla }\wedge (\vec{\nabla} \wedge \vec{\xi})= \vec{\nabla}(\vec{\nabla}.\vec{\xi})- \Delta \vec{\xi}$, we discover the existence of two types of {\it elastic waves} $\vec{A}\,exp \,i \,(\vec{k}.\vec{x}- \omega t)$ with {\it wave vector} $\vec{k}$, {\it period} $T$, {\it pulsation} $\omega=2\pi/T$ along standard notations, namely the {\it longitudinal} and {\it transversal} waves with different speeds $v_T < v_L$, which are really existing because that are responsible for earthquakes:  \\
\noindent
 \[  \left\{   \begin{array}{rcccccl}
  \vec{\nabla}.\vec{\xi}=0           &  \hspace{5mm}\Rightarrow \hspace{5mm}\vec{k}.\vec{A}=0          &\hspace{5mm}\Rightarrow &\hspace{5mm}&  \mu \triangle \vec{\xi}= \vec{f} &  \hspace{5mm} \Rightarrow \hspace{5mm} &  v_T=\sqrt{\frac{\mu}{\rho}}  \\
  \vec{\nabla}\wedge \vec{\xi}=0 & \hspace{5mm}\Rightarrow \hspace{5mm}\vec{k}\wedge \vec{A}=0 & \hspace{5mm}\Rightarrow  &\hspace{5mm}& (\lambda + 2 \mu )\triangle \vec{\xi}= \vec{f} & \hspace{5mm}   \Rightarrow \hspace{5mm} & v_L=\sqrt{\frac{\lambda + 2\mu}{\rho}}  
\end{array} \right. \] 

These comments pushed the author to a systematic use of the {\it formal adjoint} of an operator.\\

Let us explain the origin of the definition of {\it extension modules} in homological algebra by means of an elementary example. With ${\partial}_{22}\xi={\eta}^2, {\partial}_{12}\xi={\eta}^1$ for $\cal{D}$, we get  ${\partial}_1{\eta}^2-{\partial}_2{\eta}^1=\zeta$ for the CC ${\cal{D}}_1$. Then $ad({\cal{D}}_1)$ is defined by ${\mu}^2=-{\partial}_1\lambda, {\mu}^1={\partial}_2\lambda$ while $ad(\cal{D})$ is defined by $\nu={\partial}_{12}{\mu}^1+{\partial}_{22}{\mu}^2$ but the CC of $ad({\cal{D}}_1)$ are simply generated by ${\nu}'={\partial}_1{\mu}^1+{\partial}_2{\mu}^2$. Using operators, we have the two differential sequences:\\  
\[  \begin{array}{ccccl}
 \xi & \stackrel{\cal{D}}{\longrightarrow} & \eta & \stackrel{{\cal{D}}_1}{\longrightarrow} & \zeta   \\
  \nu& \stackrel{ad(\cal{D})}{\longleftarrow} & \mu & \stackrel{ad({\cal{D}}_1)}{\longleftarrow} & \lambda \\
       &    \swarrow &  &       &  \\
    \hspace*{2mm} {\nu}' &  &  &  &
  \end{array}  \]
where ${\cal{D}}_1$ generates the CC of ${\cal{D}}$ in the upper sequence but $ad({\cal{D}})$ does not generate the CC of $ad({\cal{D}}_1)$ in the lower sequence, even though ${\cal{D}}_1\circ {\cal{D}}=0 \Rightarrow ad({\cal{D}}) \circ ad({\cal{D}}_1)=0$, contrary to what happens in the Poincar\'{e} sequence for the exterior derivatve used in electromagnetism for exhibiting the Maxwell equations when $n=4$. We shall see that this comment brings the need to introduce the {\it first extension module} ${ext}^1(M)$ of the differential module $M=coker({\cal{D}})$ determined by ${\cal{D}}$.    \\

In the general case, using the same notation for a vector bundle and its set of local sections, when $ E \stackrel{ {\cal{D}}}{\longrightarrow} F$ is a given operator, its formal adjoint is ${\wedge}^nT^*\otimes E^* \stackrel{ad({\cal{D}})}{\longleftarrow} {\wedge}^nT^*\otimes F^*$ where $E^*$ and $F^*$ are respectively obtained from $E$ and $F$ by inverting the transition matrices, like $T$ and $T^*$.  \\

Before going ahead, let us prove that there may be mainly two types of differential sequences, the {\it Janet sequence} introduced by M. Janet in 1920, having to do with the tools we have studied, and a different sequence called {\it Spencer sequence} introduced by D. C. Spencer in 1970 with totally different operators ([4, 30]). For this, if $E$ is a vector bundle over the base $X$, we introduce the $q$-jet bundle $J_q(E)$ with sections ${\xi}_q: (x)\rightarrow ({\xi}^k(x), {\xi}^k_i(x), {\xi}^k_{ij}(x), ... )$ transforming like the sections $j_q(\xi):(x) \rightarrow ({\xi}^k(x), {\partial}_i{\xi}^k(x), {\partial}_{ij}{\xi}^k_(x), ...) $. The {\it Spencer operator} $D:J_{q+1}(E) \rightarrow T^*\otimes J_q(T)$ allows to compare these sections by considering the differences $({\partial}_i{\xi}^k(x)-{\xi}^k_i(x), {\partial}_i{\xi}^k_j(x) - {\xi}^k_{ij}(x), ...)$ and so on. When $\omega$ is a nondegenerate metric with Christoffel symbols $\gamma$ and Levi-Civita isomorphism $j_1(\omega)\simeq (\omega, \gamma)$ while $T=T(X)$ is the tangent bundle to $X$, we consider the second order involutive system $R_2\subset J_2(T)$ defined by considering the first order Killing system ${\cal{L}}(\xi)\omega=0$, adding its first prolongation ${\cal{L}}(\xi)\gamma=0$ and using ${\xi}_2$ instead of $j_2(\xi)$. Looking for the first order generating {\it compatibility conditions} (CC) ${\cal{D}}_1$ of the corresponding second order operator operator ${\cal{D}}$ just described, we may then look for the generating CC ${\cal{D}}_2$ of ${\cal{D}}_1$ and so on. We may proceed similarly for the injective operator $T \stackrel{j_2}{\longrightarrow} C_0(T)=J_2(T)$, finding successively $C_0(T) \stackrel{D_1}{\longrightarrow}C_1(T)$ and $C_1(T) \stackrel{D_2}{\longrightarrow}C_2(T)$ induced by $D$. When $n=2$ and $\omega$ is the $Euclidean$ metric, we have a Lie group of isometries with the $3$ infinitesimal generators $\{{\partial}_1, {\partial}_2, x^1{\partial}_2 - x^2{\partial}_1\}$. If we now consider the Weyl group defined by ${\cal{L}}(\xi)\omega = A\omega$ with $A=cst$ and ${\cal{L}}(\xi)\gamma=0$, we have to add the only dilatation $x^1{\partial}_1 + x^2 {\partial}_2$ and get the strict inclusions $R_2 \subset {\tilde{R}}_2 \subset J_2(T)$. As for the conformal system ${\hat{R}}_3 \subset J_3(T)$, according to ([25]), we have to add the two elations ${\theta}^1 = \frac{1}{2} ((x^1)^2  + (x^2)^2) {\partial}_1 + x^1 x^2 {\partial}_2$ and ${\theta}^2$ obtained by exchanging $x^1$ with $x 2$. Both systems have vanishing third symbol and we have the strict inclusions $R_3 \subset {\tilde{R}}_3   \subset {\hat{R}}_3 \subset J_3(T)$ with respective dimensions $3 < 4 < 6 $.   \\                                                                                                                                                                                                                                                                                                                                                                                                                                                                                                                                                                                                                                                                                             Collecting the results and exhibiting the induced kernel upper differential sequence, we get the following commutative {\it fundamental diagram I} where the upper down arrows are monomorphisms while the lower down arrows are epimorphisms ${\Phi}_0, {\Phi}_1, {\Phi}_2$ ([9, 12, 23]):  

 \[  \begin{array}{rccccccccccccr}
   &    &  & & & 0 & & 0 &  &0 & \\
     &   &   &    &      & \downarrow   & & \downarrow  & & \downarrow &      \\
  & 0& \longrightarrow& \hat{\Theta} &\stackrel{j_3}{\longrightarrow} & 6 &\stackrel{D_1}{\longrightarrow}& 12 &\stackrel{D_2}{\longrightarrow} &  6 &\longrightarrow  0 &  \\
  
  & 0& \longrightarrow& \tilde{\Theta} &\stackrel{j_3}{\longrightarrow}& 4 &\stackrel{D_1}{\longrightarrow}& 8 &\stackrel{D_2}{\longrightarrow} &  4 &\longrightarrow  0 &  \\

   & 0& \longrightarrow& \Theta &\stackrel{j_3}{\longrightarrow}& 3 &\stackrel{D_1}{\longrightarrow}& 6 &\stackrel{D_2}{\longrightarrow} &  3 &\longrightarrow  0 & \hspace{3mm}Spencer  \\
  &&&&& \downarrow & & \downarrow & & \downarrow & &    & \\
   & 0 & \longrightarrow &   2   & \stackrel{j_3}{\longrightarrow} &  20  & \stackrel{D_1}{\longrightarrow} &   30  &\stackrel{D_2}{\longrightarrow} &  12   &   \longrightarrow 0 &\\
   & & & \parallel && \hspace{5mm}\downarrow {\Phi}_0 & &\hspace{5mm} \downarrow {\Phi}_1 & & \hspace{5mm}\downarrow {\Phi}_2 &  &\\
   0 \longrightarrow & \Theta &\longrightarrow &   2  & \stackrel{\cal{D}}{\longrightarrow} &  17  & \stackrel{{\cal{D}}_1}{\longrightarrow} & 24 & \stackrel{{\cal{D}}_2}{\longrightarrow} &   9  & \longrightarrow  0 & \hspace{7mm} Janet \\
   
    0 \longrightarrow & \tilde{\Theta} &\longrightarrow &   2  & \stackrel{{\cal{D}}}{\longrightarrow} &  16  & \stackrel{{\cal{D}}_1}{\longrightarrow} & 22  & \stackrel{{\cal{D}}_2}{\longrightarrow} &   8  & \longrightarrow  0 & \\
    
      0 \longrightarrow & \hat{\Theta} &\longrightarrow &   2  & \stackrel{{\cal{D}}}{\longrightarrow} &  14 & \stackrel{{\cal{D}}_1}{\longrightarrow} & 18  & \stackrel{{\cal{D}}_2}{\longrightarrow} &   6  & \longrightarrow  0 & \\
     &   &   &    &      & \downarrow   & & \downarrow  & & \downarrow &      \\
       &    &  & & & 0 & & 0 &  &0 & 
   \end{array}     \]
    
\noindent
It follows that "{\it Spencer and Janet play at see-saw} ", the dimension of each {\it Janet bundle} being decreased by the same amount as the dimension of the corresponding {\it Spencer bundle} is increased. The Poincar\'{e} sequence for the exterior derivative $d$ is  ${\wedge}^0 T^* \stackrel{d}{\longrightarrow} {\wedge}^1T^* \stackrel{d}{\longrightarrow} {\wedge}^2 T^* \rightarrow 0 $ but it is only at the end of the paper that we shall understand the link with Maxwell equations when $n=4$.  \\

\noindent
{\bf 3) DIFFERENTIAL MODULES}  \\

Let $A$ be a {\it unitary ring}, that is $1,a,b\in A \Rightarrow a+b,ab \in A, 1a=a1=a$ and even an {\it integral domain} ($ab=0\Rightarrow a=0$ or $b=0$) with {\it field of fractions} $K=Q(A)$. However, we shall not always assume that $A$ is commutative, that is $ab$ may be different from $ba$ in general for $a,b\in A$. We say that $M={}_AM$ is a {\it left module} over $A$ if $x,y\in M\Rightarrow ax,x+y\in M, \forall a\in A$ or a {\it right module} $M_B$ over $B$ if the operation of $B$ on $M$ is $(x,b)\rightarrow xb, \forall b\in B$. If $M$ is a left module over $A$ and a right module over $B$ with $(ax)b=a(xb), \forall a\in A,\forall  b\in B, \forall x\in M$, then we shall say that $M={ }_AM_B$ is a {\it bimodule}. Of course, $A={ }_AA_A$ is a bimodule over itself. We define the {\it torsion submodule} $t(M)=\{x\in M\mid \exists 0\neq a\in A, ax=0\}\subseteq M$ and $M$ is a {\it torsion module} if $t(M)=M$ or a {\it torsion-free module} if $t(M)=0$. We denote by $hom_A(M,N)$ the set of morphisms $f:M\rightarrow N$ such that $f(ax)=af(x)$. We finally recall that a sequence of modules and maps is exact if the kernel of any map is equal to the image of the map preceding it. When $A$ is commutative, $hom(M,N)$ is again an $A$-module for the law $(bf)(x)=f(bx)$ as we have $(bf)(ax)=f(bax)=f(abx)=af(bx)=a(bf)(x)$. In the non-commutative case, things are more complicate and, given ${}_AM$ and ${}_AN_B$, then $hom_A(M,N)$ becomes a right module over $B$ for the law $(fb)(x)=f(x)b$ (See [10, 14, 18] for more details or [8, 14, 15, 30] for homological algebra) . \\

\noindent
{\bf DEFINITION 3.1}: A module $F$ is said to be {\it free} if it is isomorphic to a (finite) power of $A$ called the {\it rank} of $F$ over $A$ and denoted by $rk_A(F)$ while the rank $rk_A(M)$ of a module $M$ is the rank of a maximum free submodule $F\subset M$. It follows from this definition that $M/F$ is a torsion module. In the sequel we shall only consider {\it finitely presented} modules, namely {\it finitely generated} modules defined by exact sequences of the type $F_1 \stackrel{d_1}{\longrightarrow} F_0 \stackrel{p}{\longrightarrow} M\longrightarrow 0$ where $F_0$ and $F_1$ are free modules of finite ranks $m_0$ and $m_1$ often denoted by $m$ and $p$ in examples. A module $P$ is called {\it projective} if there exists a free module $F$ and another (projective) module $Q$ such that $P\oplus Q\simeq F$.\\

\noindent
{\bf PROPOSITION 3.2}: For any short exact sequence $0\rightarrow M' \stackrel{f}{\longrightarrow} M \stackrel{g}{\longrightarrow} M" \rightarrow 0$, we have the important relation $rk_A(M)=rk_A(M')+rk_A(M")$, even in the non-commutative case. As a byproduct, if $M$ admits a finite length free {\it resolution}  \,$ ... \, \stackrel{d_2}{\longrightarrow} F_1 \stackrel{d_1}{\longrightarrow} F_0 \stackrel{p}{\longrightarrow} M \rightarrow 0$, we may define the {\it Euler-Poincar\'{e} characteristic} \,${\chi}_A(M)={\sum}_r(-1)^r rk_A(F_r)=rk_A(M)$.  \\

The following classical proposition will be used later on for exhibiting the Ricci tensor and the Weyl tensor from the Riemann tensor:\\

\noindent
{\bf PROPOSITION 3.3}: We shall say that the following short exact sequence {\it splits} if one of the following equivalent three conditions holds:  \\
\[    0 \longrightarrow M' \stackrel{\stackrel{u}{\longleftarrow}}{\stackrel{f}{\longrightarrow}} M
     \stackrel{\stackrel{v}{\longleftarrow}}{\stackrel{g}{\longrightarrow}} M''   \longrightarrow 0  \]
$\bullet$ There exists a monomorphism $v:M''\rightarrow M$ called {\it lift} of $g$ and such that $g\circ v=id_{M''}$ .\\
$\bullet$ There exists an epimorphism $u:M\rightarrow M'$ called {\it lift} of $f$ and such that $u\circ f=id_{M'}$.\\
$\bullet$ There exist isomorphisms $\varphi=(u,g):M\rightarrow M'\oplus M''$ and $\psi=f+v:M'\oplus M''\rightarrow M$ that are inverse to each other and provide an isomorphism $M\simeq M'\oplus M''$ with $f\circ u+v\circ g=id_M$ and thus $ker(u)=im(v)$.  \\
These conditions are automatically satisfied if $M"$ is free or projective.  \\

Using the notation $M^*=hom_A(M,A)$, for any morphism $f:M\rightarrow N$, we shall denote by $f^*:N^*\rightarrow M^*$ the morphism which is defined by  $f^*(h)=h\circ f, \forall h\in hom_A(N,A)$ and satisfies $rk_A(f)=rk_A(im(f))=rk_A(f^*),\forall f\in hom_A(M,N)$. We may take out $M$ in order to obtain the {\it deleted sequence} $... \stackrel{d_2}{\longrightarrow} F_1 \stackrel{d_1}{\longrightarrow} F_0 \longrightarrow 0$ and apply  $hom_A(\bullet,A)$ in order to get the sequence $... \stackrel{d^*_2}{\longleftarrow} F^*_1 \stackrel{d^*_1}{\longleftarrow} F^*_0 \longleftarrow 0$. \\

\noindent
{\bf PROPOSITION 3.4}: If we define the {\it extension modules}  $ext^0_A(M)=ker(d^*_1)=hom_A(M,A)=M^*$ and $ext^i(M)=ext^i_A(M)=ker(d^*_{i+1})/im(d^*_i), \forall i\geq 1$, they do not depend on the resolution chosen and are torsion modules for $i\geq 1$. \\

We now turn to the operator framework with modules over the ring $D=K[d_1, ...,d_n]=K[d]$ of differential operators with coefficients in a differential 
field $K$ with $n$ commuting derivations $({\partial}_1,...,{\partial}_n)$, also called $D$-modules.  \\

\noindent
{\bf DEFINITION 3.5}: If a differential operator $\xi \stackrel{\cal{D}}{\longrightarrow} \eta$ is given, a {\it direct problem} is to find generating {\it compatibility conditions} (CC) as an operator $\eta \stackrel{{\cal{D}}_1}{\longrightarrow} \zeta $ such that ${\cal{D}}\xi=\eta \Rightarrow {\cal{D}}_1\eta=0$. Conversely, given $\eta \stackrel{{\cal{D}}_1}{\longrightarrow} \zeta$, the {\it inverse problem} will be to look for $\xi \stackrel{\cal{D}}{\longrightarrow} \eta$ such that ${\cal{D}}_1$ generates the CC of ${\cal{D}}$ and we shall say that ${\cal{D}}_1$ {\it is parametrized by} ${\cal{D}}$ {\it if such an operator} ${\cal{D}}$ {\it is existing}. \\

Introducing the morphism $\epsilon: M \rightarrow M^{**}$ such that $\epsilon (m)(f)=f(m), \forall m\in M, \forall f\in M^*$ and defining the differential module $N$  from $ad({\cal{D}}_1)$ exactly like we defined the differential module $M$ from ${\cal{D}}$, we finally notice that any operator is the adjoint of a certain operator because $ad(ad(P))=P, \forall P \in D$ and we get ([14, 15]):  \\

\noindent
{\bf THEOREM 3.6}: ({\it reflexivity test}) In order to check whether $M$ is {\it reflexive} or not, that is to find out a parametrization if $t(M)=0$ which {\it can be again parametrized}, the test has 5 steps which are drawn in the following diagram where $ad({\cal{D}})$ generates the CC of $ad({\cal{D}}_1)$ and ${\cal{D}}_1'$ generates the CC of ${\cal{D}}=ad(ad({\cal{D}}))$ while $ad({\cal{D}}_{-1})$ generates the CC of $ad({\cal{D}})$ and ${\cal{D}}'$ generates the CC of ${\cal{D}}_{-1}$:  \\
\[  \begin{array}{rcccccccl}
 & & & & & {\eta}'     & &  {\zeta}' &\hspace{15mm} 5  \\
 & & & &  \stackrel{{\cal{D}}'}{\nearrow}   & & \stackrel{{\cal{D}}'_1}{\nearrow} &  &  \\
4 \hspace{15mm}&\phi & \stackrel{{\cal{D}}_{-1}}{\longrightarrow}& \xi  & \stackrel{{\cal{D}}}{\longrightarrow} &  \eta & \stackrel{{\cal{D}}_1}{\longrightarrow} & \zeta &\hspace{15mm}   1  \\
 &  &  &  &  &  &  &  &  \\
 &  &  &  &  &  &  &  &  \\
 3 \hspace{15mm}& \theta &\stackrel{ad({\cal{D}}_{-1})}{\longleftarrow}& \nu & \stackrel{ad({\cal{D}})}{\longleftarrow} & \mu & \stackrel{ad({\cal{D}}_1)}{\longleftarrow} & \lambda &\hspace{15mm} 2
  \end{array}  \]
\[{\cal{D}}_1 \,\,\,parametrized \,\,\,by \,\,\,{\cal{D}} \Leftrightarrow {\cal{D}}_1={\cal{D}}'_1  \Leftrightarrow ext^1(N)=0  \Leftrightarrow  \epsilon \,\,\, injective  \Leftrightarrow t(M)=0\] 
\[{\cal{D}} \,\,\,parametrized \,\,\,by \,\,\,{\cal{D}}_{-1} \Leftrightarrow {\cal{D}}={\cal{D}}' \Leftrightarrow ext^2(N)=0 \Leftrightarrow \epsilon \,\,\, surjective  \hspace{17mm}  \]

\noindent
{\bf COROLLARY 3.7}: In the differential module framework, if $F_1 \stackrel{{\cal{D}}_1}{\longrightarrow} F_0 \stackrel{p}{\longrightarrow} M \rightarrow 0$ is a finite free presentation of $M=coker({\cal{D}}_1)$ with $t(M)=0$, then we may obtain an exact sequence $F_1 \stackrel{{\cal{D}}_1}{\longrightarrow} F_0 \stackrel{{\cal{D}}}{\longrightarrow} E $ of free differential modules where ${\cal{D}}$ is the parametrizing operator. However, there may exist other parametrizations $F_1 \stackrel{{\cal{D}}_1}{\longrightarrow} F_0 \stackrel{{\cal{D}}'}{\longrightarrow} E' $ called {\it minimal parametrizations} such that $coker({\cal{D}}')$ is a torsion module and we have thus $rk_D(M)=rk_D(E')$.  \\

\noindent
{\bf EXAMPLE 3.8}: When $n=3$, the $div $ operator can be parametrized by the $curl$ operator which can be itself parametrized by the $grad$ operator. However, using $({\xi}^1,{\xi}^2, {\xi}^3=0)$, we may obtain the new minimal parametrization $ - {\partial}_3{\xi}^2={\eta}^1, {\partial}_3{\xi}^1={\eta}^2, {\partial}_1{\xi}^2-{\partial}_2{\xi}^1={\eta}^3\Rightarrow {\partial}_1{\eta}^1+{\partial}_2{\eta}^2+{\partial}_3{\eta}^3=0$ which cannot be again parametrized.  \\

\noindent
{\bf THEOREM3.9}: the Einstein operator is self-adjoint (with a slight abuse of language), where by Einstein operator we mean the linearization of the Einstein tensor 
over the locally constant Minkowski metric $\omega$.  \\

\noindent
{\it Proof}: First of all, the linearizations of the Christoffel symbols ${\gamma}^k_{ij}$ and the Riemann tensor ${\rho}^k_{l,ij}$ are:  \\
\[  {\Gamma}^k_{ij}=\frac{1}{2} {\omega}^{kr} (d_i {\Omega}_{rj} + d_j {\Omega}_{ir} - d_r {\Omega}_{ij})  \Rightarrow  R^k_{l,ij}=d_i{ \Gamma}^k_{lj} - d_j {\Gamma}^k_{li}. \]
Setting $\Omega= {\omega}^{rs}{\Omega}_{rs}$, we deduce ${\Gamma}^r_{ri}= \frac{1}{2} d_i \Omega $ and get:  \\
\[ 2 R_{ij} = {\omega}^{rs} d_{rs} {\Omega}_{ij} + d_{ij} \Omega - {\omega}^{rs}(d_{ri}{\Omega}_{sj} + d_{rj} {\Omega}_{si})  \Rightarrow R = {\omega}^{rs} d_{rs} \Omega - {\omega}^{rs} {\omega}^{ij} d_{ri} {\Omega}_{sj}  \]
Setting $E_{ij}=R_{ij}- \frac{1}{2} {\omega}_{ij} R$ with $R ={\omega}^{rs}R_{rs} $, we obtain the linear Einstein operator (6 terms):  \\
\[  2 E_{ij} = {\omega}^{rs} d_{rs} {\Omega}_{ij} + d_{ij} \Omega - {\omega}^{rs} (d_{ri} {\Omega}_{sj} + d_{rj} {\Omega}_{si}) - 
{\omega}_{ij} ( {\omega}^{rs}d_{rs} \Omega - {\omega}^{ru} {\omega}^{sv} d_{rs} {\Omega}_{uv})  \]
It is essential to notice that the {\it Ricci operator is not self-adjoint} because we have for example:  \\
\[  {\lambda}^{ij} ({\omega}^{rs}d_{ij}{\Omega}_{rs}) \stackrel{ad}{\longrightarrow} ({\omega}^{rs}d_{ij}{\lambda}^{ij} {\Omega}_{rs} = 
({\omega}^{ij} d_{rs}{\lambda}^{rs}) {\Omega}_{ij} \]
and $ad$ provides a term appearing in $- {\omega}_{ij} R$ but {\it not} in $2R_{ij}$.  \\
After two integrations by parts, we obtain successively:  \\
\[ {\Omega}_{ij} {\omega}^{rs}  d_{rs} {\lambda}^{ij} + \Omega d_{ij} {\lambda}^{ij} - {\omega}^{rs} ({\Omega}_{sj}d_{ri} {\lambda}^{ij} + {\Omega}_{si} d_{rj} {\lambda}^{ij}) - 
     {\omega}_{ij} \Omega {\omega}^{rs} d_{rs} {\lambda}^{ij} + {\omega}_{ij} {\omega}^{ru} {\omega}^{sv} {\Omega}_{uv} d_{rs} {\lambda}^{ij} \]
Setting $\lambda = {\omega}_{ij} {\lambda}^{ij}$, we may change the indices in order to factor out ${\Omega}_{ij}$ and finally get: \\
\[  {\omega}^{rs} d_{rs} {\lambda}^{ij} + {\omega}^{ij} d_{rs} {\lambda}^{rs} - ({\omega}^{ri} d_{rs} {\lambda}^{sj} + {\omega}^{rj} d_{rs} {\lambda}^{si}) - 
{\omega}^{ij} {\omega}^{rs} d_{rs} \lambda + {\omega}^{ri} {\omega}^{sj} d_{rs} \lambda  \]
the $6$ terms being exchanged between themselves with $(1, 2, 3, 4, 5, 6) \rightarrow (1, 6, 3, 4, 5, 2)$.  \\
\hspace*{12cm}    $ \Box $ \\

\noindent
{\bf EXAMPLE 3.10}: When $n=3$ and the Euclidean metric, we obtain:  \\
\[      R = (d_{11}{\Omega}_{22} + d_{11}{\Omega}_{33} + d_{22}{\Omega}_{11} + d_{22}{\Omega}_{33} + d_{33}{\Omega}_{11} + d_{33}{\Omega}_{22}) 
    - 2 ( d_{12}{\Omega}_{12}+ d_{13}{\Omega}_{13} + d_{23}{\Omega}_{23})    \]
 \[  2 R_{12} = 2 E_{12} =  d_{33}{\Omega}_{12} + d_{12}{\Omega}_{33} - d_{13}{\Omega}_{23} - d_{23}{\Omega}_{13}  \]
\[  2 R_{11} = (d_{22} + d_{33}) {\Omega}_{11} + d_{11}({\Omega}_{22} + {\Omega}_{33}) - 2 (d_{12}{\Omega}_{12} + d_{13} {\Omega}_{13})   \]
\[  - 2 E_{11}= d_{22}{\Omega}_{33} + d_{33}{\Omega}_{22} - 2 d_{23}{\Omega}_{23}    \]
We let the reader check that eight to twelve terms are disappearing each time, a reason for which nobody saw that the Einstein equations had been written {\it exactly} (up to sign) by E. Beltrami  in 1892  in order to parametrize the Cauchy stress equations while using 6 stress functions ${\phi}_{ij} = {\phi}_{ji}$ in place of ${\Omega}_{ij} = {\Omega}_{ji}$, 25 years before Einstein ([1, 7, 23, 28]). {\it The comparison needs no comment !}.  \\

\noindent
{\bf REMARK 3.11}:  When $n=3$, we have already noticed in the first page of ([18]) that one has to take into account the factors "2" in the duality summation with Lagrange multipliers  $({\lambda}^{ij}={\lambda}^{ji})$ which is used in order to exhibit the operator $ad(Einstein)$, namely:                          \\
\[   {\lambda}^{ij} E_{ij}= {\lambda}^{11} E_{11} + 2 {\lambda}^{12} E_{12} + 2 {\lambda}^{13} E_{13} + {\lambda}^{22} E_{22} + 
2 {\lambda}^{23} E_{23} + {\lambda}^{33} E_{33} \]
In the present situation, this is sufficient in order to obtain a self-adjoint operator as follows:  \\
\[  \fbox{  $  \frac{1}{2} \left(  \begin{array}{cccccc} 
 0 & 0 & 0 & - d_{33} &  2d_{23} & - d_{22} \\
 0 &  2d_{33} &  - 2d_{23} & 0 &  - 2d_{13} &  2d_{12}  \\
 0 &  - 2d_{23} &  2d_{22} &  2d_{13} & -2d_{12} & 0 \\
 - d_{33}& 0 &  2 d_{13} & 0 & 0 & - d_{11}  \\
   2d_{23} &  - 2d_{13} &  - 2d_{12}& 0 &  2d_{11} & 0 \\
 - d_{22} &  2 d_{12} & 0 & - d_{11}& 0 & 0 
 \end{array} \right) 
 \left( \begin{array}{r}
 {\Omega}_{11} \\
 {\Omega}_{12}  \\
 {\Omega}_{13}   \\
 {\Omega}_{22}  \\
 {\Omega}_{23}  \\
 {\Omega}_{33}
 \end{array} \right) = 
 \left( \begin{array}{r}
  E_{11}  \\
  2 E_{12} \\
  2 E_{13}  \\
  E_{22}  \\
  2 E_{23}   \\
  E_{33}
 \end{array} \right)  $  }   \]    \\

We shall finally prove below that the {\it Einstein parametrization} of the stress equations is neither canonical nor minimal in the following diagrams ([23]):  \\
\[   \begin{array}{rcccccccccl}
 &4 & \stackrel{Killing}{\longrightarrow} & 10 & \stackrel{Riemann}{\longrightarrow} & 20 & \stackrel{Bianchi}{\longrightarrow} & 20 & \longrightarrow & 6 & \rightarrow 0 \\
  &   &                                                            & \parallel &       & \downarrow  &  & \downarrow &  &   \\
 &    &                                                            &  10     & \stackrel{Einstein}{\longrightarrow}  & 10  & \stackrel{div}{\longrightarrow} & 4    & \rightarrow & 0                      &   \\  
  &   &  &  &  &  &  &  &  \\   
0 \leftarrow & 4 & \stackrel{Cauchy}{\longleftarrow} & 10 & \stackrel{Beltrami}{\longleftarrow} & 20 & \stackrel{Lanczos}{\longleftarrow} & 20 &  & & \\
                     &           &                                                         & \parallel &    & \uparrow  &  &  &    \\
  &    &      & 10 &  \stackrel{Einstein}{\longleftarrow} & 10 &  &  &  & & 
\end{array}   \]

The upper $div$ induced by $Bianchi$ has {\it strictly nothing to do} with the lower $Cauchy$ operator, contrary to what is still believed today while the $10$ {\it on the right} of the lower diagram has {\it strictly nothing to do} with the perturbation of a metric which is the $10$ {\it on the left} in the upper diagram. It also follows that the Einstein equations in vacuum cannot be parametrized as we have the following diagram of operators recapitulating the five steps of the parametrizability criterion (A computer agebra exhibition of this result can be easily provided): \\
\[  \begin{array}{rcccl}
  &  &  &\stackrel{Riemann}{ }  & 20   \\
  & &  & \nearrow &    \\
 4 &  \stackrel{Killing}{\longrightarrow} & 10 & \stackrel{Einstein}{\longrightarrow} & 10  \\
  & & & &  \\
 4 & \stackrel{Cauchy}{\longleftarrow} & 10 & \stackrel{Einstein}{\longleftarrow} & 10 
\end{array}  \] \\

As a byproduct, we are facing {\it only two} possibilities, both leading to a contradiction ([27]):  \\

1) If we use the operator $S_2T^* \stackrel{Einstein}{\longrightarrow} S_2T^*$ in the "geometrical " setting of H. Poincar\'{e}, the $S_2T^*$ on the left has indeed {\it someting to do} with the perturbation ${\Omega}$ of the metric $\omega$ {\it by definition} but the $S_2T^*$ on the right has {\it strictly nothing to do} with the stress. \\

2) If we use the adjoint operator ${\wedge}^nT^*\otimes S_2T\stackrel{Einstein}{\longleftarrow} {\wedge}^nT^*\otimes S_2T$ in the "physical " setting of H. Poincar\'{e}, then ${\wedge}^nT^*\otimes S_2T$ on the left has of course {\it something to do} with the stress ${\sigma}$ but the ${\wedge}^nT^*\otimes S_2T$ on the right has {\it strictly nothing to do} with a perturbation of the metric as it is just the definition of the stress functions ${\phi}$ allowing to parametrize the Cauchy operator. \\

{\it These purely mathematical results question the origin and existence of gravitational waves}.  \\

We shall finally prove below, from what has been said in the beginning of this paper, that the parametrization of the $Cauchy$ operator by $ad(Ricci)$ is neither canonical nor minimal in the following diagrams:  \\
\[   \begin{array}{rcccccccccl}
 &4 & \stackrel{Killing}{\longrightarrow} & 10 & \stackrel{Riemann}{\longrightarrow} & 20        \\                                                      
  &      &                    &  \parallel &   &\downarrow  &  &   \\
 &    &                                                            &  10     & \stackrel{Ricci}{\longrightarrow}  & 10      \\   
  &  &  &  &  &  \\
0 \leftarrow & 4 & \stackrel{Cauchy}{\longleftarrow} & 10 & \stackrel{Beltrami}{\longleftarrow} & 20  \\
                     &           &                                                         & \parallel &    & \uparrow     \\
  &    &      & 10 &  \stackrel{ad(Ricci)}{\longleftarrow} & 10  
\end{array}   \]

\noindent
 {\bf COROLLARY 3.12}: When $n\geq 4$, as the extension modules are torsion modules, each component of the Weyl tensor is a torsion element killed by the $Dalembert$ operator whenever the Einstein equations in vacuum are satisfied by the metric. Hence, there exists a second order  operator ${\cal{Q}}$ such that we have an identity describing the so-called {\it Lichnerowicz wave equations} (in France !):  \\
 \[   \Box \circ Weyl = {\cal{Q}} \circ Ricci    \]

 \noindent
{\bf 4) GAUGE THEORY}   \\

\noindent
{\it Gauging procedure} : If $y=a(t)x+b(t)$ with $a(t)$ a time depending orthogonal matrix ({\it rotation}) and $b(t)$ a time depending vector ({\it translation}) describes the movement of a rigid body in ${\mathbb{R}}^3$, then the projection of the speed $v=\dot{a}(t) x + \dot{b}(t)$ in an orthogonal frame fixed in the body is 
$a^{-1}v=(a^{-1}\dot{a}) x + a^{-1}\dot{b}$ and the kinetic energy is a quadratic function of the $1$-forms $a^{-1}\dot{a}$ and $a^{-1}\dot{b}$.  \\

More generally, we may consider a map $a: X \longrightarrow G:x \longrightarrow a(x)$, introduce the tangent mapping $T(a):T=T(X) \longrightarrow T(G):dx \longrightarrow da=\frac{\partial a}{\partial x} dx$ and consider the family of left invariant $1$-forms $a^{-1}da=A=(A^{\tau}_i(x)dx^i)$ with value in the {\it Lie algebra} ${\cal{G}}=T_e(G)$, the tangent space of $G$ at the identity $e\in G$ with {\it structure constants} $c=(c^{\tau}_{\rho\sigma})$ for $ 1 \leq \rho, \sigma, \tau \leq p $. We may introduce the $2$-forms ${\partial}_iA^{\tau}_j-{\partial}_jA^{\tau}_i-c^{\tau}_{\rho\sigma}A^{\rho}_iA^{\sigma}_j=F^{\tau}_{ij}$ with value in ${\cal{G}}$, simply denoted by $dA-[A,A]=F$, and we have $A=a^{-1}da \Leftrightarrow F=0$ by pulling back on $X$ the {\it Maurer-Cartan equations} (MC) on $G$ ([11]). \\

In 1956, at the birth of GT, the above notations were coming from the EM potential $A$ and EM field $dA=F$ of relativistic Maxwell theory. 
Indeed, $G=U(1)$ (unit circle in the complex plane)$\longrightarrow dim ({\cal{G}})=1$ was the {\it only possibility} to get a $1$-form $A$ and a $2$-form $F$ when 
$c=0$. Such a choice convinced people that the EM field F should be related to the curvature introduced by Cartan, as a $2$-form with value in a Lie algebra, exactly like the Riemann tensor is considered as a $2$-form with value in rotation and that {\it curvature + torsion} should be therefore a generalization of curvature alone, roughly that the "{\it field}" should be a section of the Spencer bundle $C_2$. \\

On the contrary, in the conformal framework, that is when $G$ {\it is acting on} $X$, the second order jets ({\it elations}) ${\xi}^k_{ij}={\delta}^k_ia_j+{\delta}^k_ja_i-{\omega}_{ij}{\omega}^{kr}a_r \Rightarrow {\xi}^r_{ri}=na_i$ behave like the $1$-form $a_i(x)dx^i$ and the corresponding part of the Spencer operator $D$ is a $1$-form with value in $1$-form, that is a $(1,1)$-covariant tensor providing the EM field as a $2$-form by skewsymmetrization. This result, namely to construct lagrangians on the image of the induced Spencer operator $D_1$, is thus {\it perfectly coherent} with rigid body dynamics, Cosserat elasticity and Maxwell theory but in {\it total contradiction} with GT because $U(1)$ is not acting on space-time and there is a {\it shift by one step} in the interpretation of the Poincar\'{e} sequence involved because the fields are now described by $1$-forms ([11, 18]).  \\

\noindent
{\it Gauging procedure revisited} : Finally, we may extend the action $y=f(x,a)$ to $y_q=j_q(f)(x,a)$ in order to eliminate the parameters when $q$ is large enough. In this case, we may set $f(x)=f(x,a(x))$ and $f_q(x)=j_q(f)(x,a(x))$ in order to obtain $a(x)=a=cst \Leftrightarrow f_q=j_q(f)$ because 
$Df_{q+1} = j_1(f_q)-f_{q+1} = (\frac{\partial f_q(x,a(x))}{\partial a^{\tau}}{\partial}_ia^{\tau}(x)) $ and the matrix involved has maximum rank $p$.  \\. \\

\noindent
{\bf 5) COSSERAT VERSUS WEYL}  \\

Computing the formal adjoint $ad(D_1)$ of the first Spencer operator $D_1$ induced by $D$ for the group of rigid motions when $n=2$, we get ([16]):  \\
\[  {\sigma}^{11}({\partial}_1{\xi}_1 - {\xi}_{1,1}) + {\sigma}^{12} ({\partial}_2{\xi}_1-{\xi}_{1,2})+{\sigma}^{21}({\partial}_1{\xi}_2-{\xi}_{2,1}) + {\sigma}^{22}({\partial}_2{\xi}_2-{\xi}_{2,2})+{\mu}^{12,r}({\partial}_r{\xi}_{1,2}-{\xi}_{1,2r})  \]
Integrating by parts with ${\xi}_{1,1}=0, {\xi}_{1,2}+{\xi}_{2,1}=0, {\xi}_{2,2}=0 \Rightarrow {\xi}_{1,2r}=0$, we obtain {\it at once} and {\it exactly} the {\it Cosserat couple-stress equations}:  \\
\[ {\partial}_1{\sigma}^{11}+{\partial}_2{\sigma}^{12}=f^1, \,\, \,{\partial}_1{\sigma}^{21}+{\partial}_2{\sigma}^{22}=f^2, \,\,\,
     {\partial}_r{\mu}^{12,r} +{\sigma}^{12} - {\sigma}^{21}=m^{12}  \]
allowing to have now a non-symmetrical stress and a new {\it first order parametrization}:  \\
\[  {\sigma}^{11}={\partial}_2{\phi}^1, \,\,{\sigma}^{12}= - {\partial}_1{\phi}^1,\,\, {\sigma}^{21}= - {\partial}_2{\phi}^2, \,\,{\sigma}^{22}= {\partial}_1{\phi}^2, \,\,  {\mu}^{12,1}={\partial}_2{\phi}^3 + {\phi}^1,  \,\,{\mu}^{12,2}= - {\partial}_1{\phi}^3 - {\phi}^2  \]
with a potential $({\phi}^1,{\phi}^2, {\phi}^3)\in {\wedge}^2T^*\otimes {\wedge}^2T\otimes R_2^*\simeq R_2^*$ with $dim(R_2)=3$. These equations can be extended by adding the only dilatation with infinitesimal generator $x^i{\partial}_i$ in order to provide the {\it virial equations} (See [11] for more details).  \\

Similarly, going along the idea pioneered by Weyl in 1918 ([32]), we obtain with $i<j$: \\
\[ ({\partial}_i {\xi}^r_{rj} - {\xi}^r_{rji}) - ({\partial}_j {\xi}^r_{ri} - {\xi}^r_{rij}) = {\partial}_i {\xi}^r_{rj} - {\partial}_j {\xi}^r_{ri} \,\,\, \Rightarrow \,\,\, {\cal{J}}^i ( {\partial}_i {\xi}^r_r - 
{\xi}^r_{ri}) + {\cal{F}}^{ij} ( {\partial}_i {\xi}^r_{rj} - {\partial}_j {\xi}^r_{ri} )  \]
An integration by parts brings the equations $ {\partial}_r {\cal{F}}^{ir} - {\cal{J}}^i = 0$ and ${\partial}_i {\cal{J}}^i=0$, a result leading to: \\

\noindent
 {\it There is no conceptual difference at all between the Cosserat couple-stress equations and the second set of Maxwell equations. Only the groups are different}. \\     
      
As a next crucial step, let us consider the Lie group of transformations of $X$ described by the action of a Lie group $G$ with local coordinates $(a^{\tau})$, identity $e\in G$ and Lie algebra ${\cal{G}}=T_e(G)$, on $X$ with infinitesimal generators ${\theta}_{\tau}={\theta}^k_{\tau}(x){\partial}_k$ and introduce the section ${\xi}_q={\lambda}^{\tau}(x)j_q({\theta}_{\tau})$ with $\lambda \in {\wedge}^0T^*\otimes {\cal{G}}$. We have thus 
$ ({\xi}^k(x)={\lambda}^{\tau}(x){\theta}^k_{\tau}(x), {\xi}^k_i(x)={\lambda}^{\tau}(x){\partial}_i{\theta}^k_{\tau}(x), ...)$ and get at once $D_1{\xi}_q=D{\xi}_{q+1}=(d{\lambda}^{\tau})j_q({\theta}_{\tau})$ where $d:{\wedge}^0T^*\otimes {\cal{G}} \rightarrow {\wedge}^1T^*\otimes {\cal{G}}$ is the exterior derivative, a result proving that the Spencer sequence is (locally) isomorphic to the tensor product by ${\cal{G}}$ of the Poincar\'{e} sequence for $d$, in a coherent way with the second Example of the Introduction. As the extension modules of a module $M$ do not depend on the resolution of $M$, it follows that, if ${\cal{D}}_{r+1}$ generates the CC of ${\cal{D}}_r$  in a Janet sequence like in the Introduction, then $ad({\cal{D}}_r)$ generates the CC of $ad({\cal{D}}_{r+1})$ while, if $D_{r+1}$ generates the CC of $D_r$ in the corresponding Spencer sequence, then $ad(D_r)$ generates the CC of $ad(D_{r+1})$, {\it though all these operators are quite different}, a result not evident at all that Lanczos and followers could have not even been able to imagine ([24, 28]).\\

It remains to prove that, in this new framework, the Ricci tensor only depends on the symbol ${\hat{g}}_2\simeq T^*\subset S_2T^*\otimes T$ of the first prolongation ${\hat{R}}_2\subset J_2(T) $ of the conformal Killing system ${\hat{R}}_1\subset J_1(T)$ with symbol ${\hat{g}}_1\subset T^*\otimes T$ defined by the equations ${\omega}_{rj}{\xi}^r_i + {\omega}_{ir}{\xi}^r_j - \frac{2}{n}{\omega}_{ij}{\xi}^r_r=0$ {\it not depending on any conformal factor}. The next commutative diagram covers both situations,  taking into account that the equations of both the classical and conformal Killing operator are homogeneous. The purely algebraic {\it Spencer map} 
$\delta: g_{q+1} \rightarrow T^* \otimes  g_q$ with symbol $g_q = R_q \cap S_qT^* \otimes E  \subset J_q(E)$ is induced by $-D$ and all the sequences are exact by definition but perhaps the left column:  \\
\[  \begin{array}{rcccccccl}
   &  0 & & 0 & & 0 &  &  &   \\
   & \downarrow & & \downarrow & & \downarrow & & &  \\
0\rightarrow & g_3 & \rightarrow &  S_3T^*\otimes T & \rightarrow & S_2T^*\otimes F_0& \rightarrow & F_1 & \rightarrow 0  \\
   & \hspace{2mm}\downarrow  \delta  & & \hspace{2mm}\downarrow \delta & &\hspace{2mm} \downarrow \delta & & &  \\
0\rightarrow& T^*\otimes g_2&\rightarrow &T^*\otimes S_2T^*\otimes T & \rightarrow &T^*\otimes T^*\otimes F_0 &\rightarrow & 0 &  \\
   &\hspace{2mm} \downarrow \delta &  &\hspace{2mm} \downarrow \delta & &\hspace{2mm}\downarrow \delta &  &  &   \\
0\rightarrow & {\wedge}^2T^*\otimes g_1 & \rightarrow & \underline{{\wedge}^2T^*\otimes T^*\otimes T} & \rightarrow & {\wedge}^2T^*\otimes F_0 & \rightarrow & 0 &  \\
   &\hspace{2mm}\downarrow \delta  &  & \hspace{2mm} \downarrow \delta  &  & \downarrow  & &  &  \\
0\rightarrow & {\wedge}^3T^*\otimes T & =  & {\wedge}^3T^*\otimes T  &\rightarrow   & 0  &  &  &   \\
    &  \downarrow  &  &  \downarrow  &  &  &  &  &  \\
    &  0  &   & 0  & &  &  &  &
\end{array}  \]

\noindent
{\bf THEOREM 5.1}: Introducing the $\delta$-cohomology bundles $H^2(g_1)$ at ${\wedge}^2T^*\otimes g_1$ and $H^2({\hat{g}}_1)$ at ${\wedge}^2T^*\otimes {\hat{g}}_1$ while taking into account that $g_1\subset {\hat{g}}_1$, $g_2=0$ ${\hat{g}}_2 \simeq T^*$ and ${\hat{g}}^3=0$, we have the commutative and exact "{\it fundamental diagram II} ":  \\
 \[ \begin{array}{rcccccccl}
 & & & & & & & 0 & \\
 & & & & & & & \downarrow & \\
  & & & & & 0& & S_2T^* &  \\
  & & & & & \downarrow & & \downarrow  &  \\
   & & & 0 &\longrightarrow & Z^2(g_1) & \longrightarrow & H^2(g_1)  & \longrightarrow 0  \\
   & & & \downarrow & & \downarrow & &  \downarrow  &  \\
   & 0 &\longrightarrow & T^*\otimes {\hat{g}}_2 & \stackrel{\delta}{\longrightarrow} & Z^2({\hat{g}}_1) & \longrightarrow & H^2({\hat{g}}_1) & \longrightarrow 0  \\
    & & & \downarrow & & \downarrow & & \downarrow     &   \\
 0 \longrightarrow & S_2T^* & \stackrel{\delta}{\longrightarrow}& T^*\otimes T^* &\stackrel{\delta}{\longrightarrow} & {\wedge}^2T^* & \longrightarrow & 0 &   \\
   & & & \downarrow &  & \downarrow & & &  \\
   & & & 0 & & 0 & & &  \\
  & & & & &  & & &     
   \end{array}  \]
The splitting sequence $ 0 \rightarrow S_2T^* \rightarrow F_1 \rightarrow {\hat{F}}_1 \rightarrow 0  $ provides a totally unusual interpretation of the successive Ricci, Riemann and Weyl tensors. It follows that $ dim(F_1) = n^2(n^2-1)/12$, \,\, $dim({\hat{F}}_1) = n(n+1)(n+2)(n-3)/12$ and the $Weyl$ operator is of order $3$ when $n=3$ but of order $2$ for $n\geq 4$. Similar results could be obtained for the $Bianchi$-type operator which is of order 4 when n=4 (See [21-24] for more details). \\

\noindent
{\bf 6) JANET VERSUS SPENCER}    \\

Whenever $R_q\subseteq J_q(E)$ is an involutive system of order $q$ on $E$, we may define the {\it Janet bundles} $F_r$ for $r=0,1,...,n$ by the short exact sequences ([9, 12]):  \\
\[0 \rightarrow {\wedge}^rT^*\otimes R_q+\delta ({\wedge}^{r-1}T^*\otimes S_{q+1}T^*\otimes E)\rightarrow {\wedge}^rT^*\otimes J_q(E) \rightarrow F_r \rightarrow 0  \]
We may pick up a section of $F_r$, lift it up to a section of ${\wedge}^rT^*\otimes J_q(E)$ that we may lift up to a section of ${\wedge}^rT^*\otimes J_{q+1}(E)$ and apply $D$ in order to get a section of ${\wedge}^{r+1}T^*\otimes J_q(E) $ that we may project onto a section of $F_{r+1}$ in order to construct an operator ${\cal{D}}_{r+1}:F_r\rightarrow F_{r+1}$ generating the CC of ${\cal{D}}_r$ in the canonical {\it linear Janet sequence}:  \\
\[  0 \longrightarrow  \Theta \longrightarrow E \stackrel{\cal{D}}{\longrightarrow} F_0 \stackrel{{\cal{D}}_1}{\longrightarrow}F_1 \stackrel{{\cal{D}}_2}{\longrightarrow} ... \stackrel{{\cal{D}}_n}{\longrightarrow} F_n \longrightarrow 0   \]
If we have two involutive systems $R_q \subset {\hat{R}}_q \subset J_q(E)$, {\it the Janet sequence for} $R_q$ {\it projects onto the Janet sequence for} ${\hat{R}}_q$ and we may define inductively {\it canonical epimorphisms} $F_r \rightarrow {\hat{F}}_r \rightarrow 0$ for 
$r=0, 1,...,n$ by comparing the previous sequences for $R_q$ and ${\hat{R}}_q$.  \\ 
A similar procedure can also be obtained if we define the Spencer bundles $C_r$ for $r=0,1,...,n$ by the short exact sequences ([9, 12]):  \\
\[ 0 \rightarrow  \delta ({\wedge}^{r-1}T^*\otimes g_{q+1} )  \rightarrow  {\wedge}^rT^*\otimes R_q  \rightarrow  C_r  \rightarrow 0 \]
We may pick up a section of $C_r$, lift it to a section of ${\wedge}^rT^*\otimes R_q$, lift it up to a section of ${\wedge}^rT^*\otimes R_{q+1}$ and apply $D$ in order to construct a section of ${\wedge}^{r+1}\otimes R_q$ that we may project to $C_{r+1}$ in order to construct an operator $D_{r+1}:C_r \rightarrow C_{r+1}$ generating the CC of $D_r$ in the canonical {\it linear Spencer sequence} which is {\it another completely different resolution} of the set $\Theta$ of (formal) solutions of $R_q$:  \\
\[    0 \longrightarrow \Theta \stackrel{j_q}{\longrightarrow} C_0 \stackrel{D_1}{\longrightarrow} C_1 \stackrel{D_2}{\longrightarrow} C_2 \stackrel{D_3}{\longrightarrow} ... \stackrel{D_n}{\longrightarrow} C_n\longrightarrow 0  \]
However, if we have two systems as above, {\it the Spencer sequence for} $R_q$ {\it is now contained into the Spencer sequence for} 
${\hat{R}}_q$ and we may construct inductively {\it canonical monomorphisms} $0\rightarrow C_r \rightarrow {\hat{C}}_r$ for $r=0,1,...,n$ by comparing the previous sequences for $R_q$ and ${\hat{R}}_q$.   \\
When dealing with applications, we have set $E=T$ and considered systems of finite type Lie equations determined by Lie groups of transformations and $ad({\cal{D}}_r)$ generates the CC of $ad({\cal{D}}_{r+1})$ while $ad(D_r)$ generates the CC of $ad(D_{r+1})$. 
We have obtained in particular $C_r={\wedge}^rT^*\otimes R_q \subset {\wedge}^rT^*\otimes {\hat{R}}_q ={\hat{C}}_r$ when comparing the classical and conformal Killing systems, but {\it these bundles have never been used in physics}. Therefore, instead of the classical Killing system $R_2\subset J_2(T)$ defined by $\Omega \equiv {\cal{L}}(\xi)\omega=0$ {\it and} $\Gamma\equiv {\cal{L}}(\xi)\gamma=0$ or the conformal Killing system ${\hat{R}}_2\subset J_2(T)$ defined by $\Omega\equiv {\cal{L}}(\xi)\omega=A(x)\omega$ and ${\Gamma} \equiv {\cal{L}}(\xi)\gamma= ({\delta}^k_iA_j(x) +{\delta} ^k_j A_i(x) -{\omega}_{ij}{\omega}^{ks}A_s(x)) \in S_2T^*\otimes T$, we may introduce the {\it intermediate differential system} ${\tilde{R}}_2 \subset J_2(T)$ defined by ${\cal{L}}(\xi)\omega=A\omega$ with $A=cst$ and $\Gamma \equiv {\cal{L}}(\xi)\gamma=0 $, for the {\it Weyl group} obtained by adding the only dilatation with infinitesimal generator $x^i{\partial}_i$ to the Poincar\'e group. We have $R_1\subset {\tilde{R}}_1={\hat{R}}_1$ but the strict inclusions $R_2 \subset {\tilde{R}}_2 \subset {\hat{R}}_2$ and we discover {\it exactly} the group scheme used through this paper, both with the need to {\it shift by one step to the left} the physical interpretation of the various differential sequences used. Indeed, as the symbol ${\hat{g}}_2\simeq T^*$ is defined by the $n(n -1)(n+2)/2$ linear equations:  \\
\[     {\xi}^k_{ij} - \frac{1}{n}( {\delta}^k_i {\xi}^r_{rj} + {\delta}^k_j {\xi}^r_{ri} - {\omega}_{ij} {\omega}^{ks} {\xi}^r_{rs}) = 0   \]
 that do not depend on any conformal factor, the first Spencer operator ${\hat{R}}_2\stackrel{D_1}{\longrightarrow} T^*\otimes {\hat{R}}_2$ is induced by the usual Spencer operator ${\hat{R}}_3 \stackrel{D}{\longrightarrow} T^*\otimes {\hat{R}}_2:(0,0,{\xi}^r_{rj},{\xi}^r_{rij}=0) \rightarrow (0,{\partial}_i0-{\xi}^r_{ri}, {\partial}_i{\xi}^r_{rj}- 0)$ and thus projects by cokernel onto the induced operator $T^* \rightarrow T^*\otimes T^*$. Composing with $\delta$, it projects therefore onto $T^*\stackrel{d}{\rightarrow} {\wedge}^2T^*:A \rightarrow dA=F$ as in EM and so on by using the fact that $D_1$ 
{\it and} $d$ {\it are both involutive}, or the composition of epimorphisms:  \\
\[ {\hat{C}}_r \rightarrow {\hat{C}}_r/{\tilde{C}}_r\simeq {\wedge}^rT^*\otimes ({\hat{R}}_2/{\tilde{R}}_2) \simeq {\wedge}^rT^*\otimes {\hat{g}}_2\simeq {\wedge}^rT^*\otimes T^*\stackrel{\delta}{\longrightarrow}{\wedge}^{r+1}T^* \]
The main result we have obtained is thus to be able to increase the order and dimension of the underlying jet bundles and groups, proving therefore that any $1$-form with value in the second order jets ${\hat{g}}_2$ ({\it elations}) of the conformal Killing system (conformal group) can be decomposed uniquely into the direct sum $(R,F)$ where $R$ is a section of the {\it Ricci bundle} $S_2T^*$ and the EM field $F$ is a section of ${\wedge}^2T^*$. This was {\it exactly} the dream of Weyl (32]). \\   

\noindent
{\bf 7) CONCLUSION}  \\

These new unavoidable methods based on the formal theory of systems of partial differential equations and Lie pseudogroups provide the common secret of the three famous books ([2, 6, 32]) published about at the same time at the beginning of the last century. Indeed, the Spencer operator can always be exhibited even if there is no group background and, when only constant sections are considered, one recovers exactly (up to sign) the operator introduced by Macaulay for studying {\it inverse systems} ([17]). As a main consequence, we have:  \\
\noindent
{\it The mathematical structures of electromagnetism and gravitation only depend on second order jets}.\\
This paper can also be considered as an elementary summary of certain recent results presented in the references below. A much more difficult non-linear version of the preceding results can be found in ([10-12, 18, 22, 26]). Finally, we hope to have convinced the reader that most of the applications presented are providing explicit examples that are tricky enough in order to justify the use of computer algebra in a near future.    \\

\vspace{1cm}
\noindent
{\bf REFERENCES}:  \\

\noindent 
[1] Beltrami, E.: Osserazioni sulla nota precedente. Atti della Accademia Nazionale dei Lincei Rend. 1, 5, 141-142 (1892) \\ 
\noindent
[2] Cosserat, E., Cosserat, F.: Th\'{e}orie des Corps D\'{e}formables. Hermann, Paris (1909)    \\
\noindent
[3] Foster, J., Nightingale, J.D.: A Short Course in General Relativity. Longman, New York (1979)  \\
 \noindent
[4] Janet, M.: Sur les syst\`{e}mes aux d\'{e}riv\'{e}es partielles. Journal de Math\'{e}matique, 8, 65-151(1920)     \\
\noindent 
[5] Kumpera, A., Spencer, D.C.: Lie Equations. Ann. Math. Studies 73, Princeton University Press, Princeton (1972).\\
 \noindent 
[6] Macaulay, F.S.: The Algebraic Theory of Modular Systems. Cambridge University Press, Cambridge (1916)  \\
\noindent
[7] Morera, G.: Soluzione generale della equazioni indefinite di equilibrio di un corpo continuo. Atti della Accademia Nazionale dei Lincei Rend. 1, 5, 137-141+233-234 (1892)   \\
\noindent
[8] Northcott, D.G.: An Introduction to Homological Algebra, Cambridge university Press, Cambridge (1966).  \\
\noindent
[9] Pommaret, J.-F.: Systems of Partial Differential Equations and Lie Pseudogroups. Gordon and Breach, New York (1978). Russian translation: MIR, Moscow (1983).\\
\noindent
[10] Pommaret, J.-F.: Differential Galois Theory. Gordon and Breach, New York (1983). \\
 \noindent
[11] Pommaret, J.-F.: Lie Pseudogroups and Mechanics. Gordon and Breach, New York (1988)     \\
\noindent
[12] Pommaret, J.-F.: Partial Differential Equations and Group Theory, Kluwer, Dordrecht (1994). https://doi.org/10.1007/978-94-017-2539-2    \\
\noindent
[13] Pommaret, J.-F.: Dualit\'{e} diff\'{e}rentielle et applications. Comptes Rendus Acad\'{e}mie des Sciences Paris, S\'{e}rie I, 320, 1225-1230 (1995)   \\
 \noindent
[14] Pommaret, J.-F.: Partial Differential Control Theory. Kluwer, Dordrecht (2001)   \\
\noindent
[15] Pommaret, J.-F.: Algebraic analysis of control systems defined by partial differential equations. In "Advanced Topics in Control Systems Theory", Springer, Lecture Notes in Control and Information Sciences, LNCIS 311, 5, 155-223 (2005)  \\
\noindent
[16] Pommaret, J.-F.: Parametrization of Cosserat Equations. Acta Mechanica 215, 43-55 (2010)   \\
 \noindent
[17] Pommaret, J.-F.: Macaulay Inverse systems revisited. Journal of Symbolic Computations 46, 1049-1069 (2011)   \\
\noindent
[18] Pommaret, J.-F.: Spencer Operator and Applications: From Continuum Mechanics to Mathematical Physics, in "Continuum Mechanics-Progress in Fundamentals and Engineering Applications", Dr. Yong Gan (Ed.), London, ISBN: 978-953-51-0447--6, InTech, 2012, Available from: \\
https://doi.org/10.5772/intechopen.35607       \\
\noindent
[19] Pommaret, J.-F.: The Mathematical Foundations of General Relativity Revisited. Journal of Modern Physics, 4, 223-239 (2013). https://doi.org/10.4236/jmp.2013.48A022   \\
\noindent
[20] Pommaret, J.-F.: Airy, Beltrami, Maxwell, Einstein and Lanczos potentials revisited. Journal of Modern Physics, 7 699-728 (2016). https://doi.org/10.4236/jmp.2016.77068   \\
\noindent
[21] Pommaret, J.-F.: Deformation Theory of Algebraic and Geometric Structures, Lambert Academic Publisher (LAP), Saarbrucken, Germany (2016). A short summary can be found in "Topics in Invariant Theory ", S\'{e}minaire P. Dubreil/M.-P. Malliavin, Springer Lecture Notes in Mathematics, 1478, 244-254 (1990).
https://arxiv.org/abs/1207.1964  \\
\noindent
[22] Pommaret, J.-F.: New Mathematical Methods for Physics. Mathematical Physics Books, Nova Science Publishers, New York (2018) \\
\noindent
[23] Pommaret, J.-F.: The mathematical foundations of elasticity and electromagnetism revisited. Journal of Modern Physics, 10 1566-1595 (2021). https://doi.org/10.4236/jmp.1013104     \\
\noindent
[24] Pommaret, J.-F.: Homological solution of the Lanczos problems in arbitrary dimension. Journal of Modern Physics, 12 , 829-858 (2021). https://doi.org/10.4236/jmp.2021.126053  \\
\noindent
[25] Pommaret, J.-F.: The conformal group revisited. https://arxiv.org/abs/2006.03449 \\
https://doi.org/10.4236/jmp.2021.1213106.  \\
\noindent
[26] Pommaret, J.-F.: Nonlinear conformal electromagnetism. https://arxiv.org/abs/2007.01710.  \\
https://doi.org/10.4236/jmp.2022.134031   \\
\noindent
[27] Pommaret, J.-F.: How Many Structure Constants do Exist in Riemannian Geometry ?. Mathematics in Computer Science 16, 23 (2022). https://doi.org/10.1007/s11786-022-00546-3      \\
\noindent
[28] Pommaret, J.-F.: Gravitational waves and parametrizations of linear differential operators. In Gravitational Waves: Theory and Observations, IntechOpen (2023)   \\
https://www.intechopen.com/online-first/1119249. https://doi.org/10.5992/intechopen.100651      \\
\noindent
[29] Pommaret, J.-F.: Killing Operator for the Kerr Metric. Journal of Modern Physics 14 31-59 (2023). https://doi.org/10.4236/jmp.2023.141003      \\
\noindent
[30] Rotman, J.J.: An Introduction to Homological Algebra. Academic Press (1979), Springer, New york (2009).   \\
\noindent
[31] Spencer, D.C.: Overdetermined systems of partial differential equations. Bulletin of the American Mathematical Society, 75, 1-114 (1965)      \\
\noindent
[32] Weyl, H.: Space, Time, Matter. Berlin (1918), Dover, New York (1952)   \\ 
\noindent
[33] Zerz, E.: Topics in Multidimensional Linear Systems Theory. In: Lecture Notes in Control and Information Sciences, LNCIS 256, Springer (2000)   \\

\end{document}